\documentclass[11pt]{amsart}
\usepackage{amssymb}

\theoremstyle{plain}
\newtheorem{theorem}{Theorem}

\newtheorem{lemma}[theorem]{Lemma}

\theoremstyle{remark}
\newtheorem{rem}{Remark}
\newtheorem*{ack}{Acknowledgement}

\newcommand{\R}{\mathbb{R}}

\newcommand{\N}{\mathbb{N}}
\newcommand{\cF}{\mathcal{F}}
\newcommand{\cG}{\mathcal{G}}
\newcommand{\eps}{\varepsilon}
\newcommand{\SL}{\mathrm{SL}(d,\mathbb{R})}
\newcommand{\GL}{\mathrm{GL}(d,\mathbb{R})}
\newcommand{\GIC}{\mathcal{G}_{\mathrm{IC}}}
\newcommand{\one}{\mathcal{G}_{\mathrm{OPS}}}
\newcommand{\dmu}{}
\newcommand{\wedgek}{\mathord{\wedge}^k}
\newcommand{\hlak}{\hat{\lambda}_k}
\DeclareMathOperator*{\esssup}{ess\,sup}

\begin{document}

\title{$L^p$-generic cocycles have one-point Lyapunov spectrum}
\date{November 1, 2002}
\author[A. Arbieto, J. Bochi]{Alexander Arbieto and Jairo Bochi}
\thanks{A. A. is supported by CNPq and J. B. is supported by CNPq-Profix}

\begin{abstract}
We show the sum of the first $k$ Lyapunov exponents of linear cocycles 
is an upper semicontinuous function  
in the $L^p$ topologies, for any $1 \le p \le \infty$ and $k$.
This fact, together with a result from Arnold and Cong, implies that
the Lyapunov exponents of the $L^p$-generic cocycle, $p<\infty$, are all equal.
\end{abstract}

\maketitle

\section{Introduction}

The Lyapunov exponents of products of random matrices,
has been the focus of much study since the pioneering work of
Furstenberg and Kesten~\cite{FurKes}.
A major problem in the theory is to describe the Lyapunov spectrum of ``typical'' cocycles.
In the case of products of independent identically distributed matrices, for example,
for ``most'' choices of the probability distribution,
all Lyapunov exponents are different, see e.g. \cite[p.~78]{BouLac}
or \cite[theorem~6.9]{GolMar}.

\smallskip

Here we consider the general setting of linear cocycles.
Fix a probability space $(X, \cF, \mu)$ and an 
ergodic invertible measure preserving transformation $T:X \to X$.
Take a measurable map $A: X \to \GL$ (satisfying an integrability condition).
Then, by the theorem of Oseledets, the Lyapunov exponents of $A$ (with multiplicity) 
$$
\lambda_1(A) \ge \ldots \ge \lambda_d(A)
$$  
are defined: 
they are the possible values of the limits
$$
\lim_{n \to \pm \infty}\frac{1}{n} \log\|A(T^{n-1} x) \cdots A(Tx) A(x) v\|,
$$
for $\mu$-almost every $x \in X$ and $v \in \R^d \smallsetminus \{0\}$.

\smallskip

One may then consider some class of cocycles $A$
and ask ``how large'' is the subset
of those cocycles with simple spectrum (i.e., with all Lyapunov exponents different),
or with one-point spectrum (all exponents equal).
By a ``large'' set, we usually mean a residual (dense $G_\delta$) set, in some topology.
Of course, the answer to the vague question above depends on
the class of the cocycles considered and of the topology.
For instance, Knill~\cite{Kni} proved that 
in the space of all bounded measurable cocycles with values in $\mathrm{SL}(2,\R)$,
the set of cocycles with two different Lyapunov exponents  
is dense in the uniform topology.
This set contains the open subset of uniformly hyperbolic cocycles.
It was showed by Bochi~\cite{Boc} that the set of cocycles that 
either are uniformly hyperbolic or have both exponents equal to zero
is a residual subset.
This result also holds for continuous cocycles,
while it is an open question whether Knill's theorem does.

Knill's result was extended (through different techniques)
to cocycles with values in $\GL$, for any $d\ge 2$,
by Arnold and Cong~\cite{ArCon2}.
Also, Bochi and Viana generalized Bochi's result to 
very general matrix groups in~\cite{BocVi}.

\smallskip
 
In this paper we consider the space of $\GL$-valued cocycles,
endowed with the weaker topologies $L^p$, $1 \le p < \infty$
(see section~\ref{s.2}).
In this case, the set of cocycles with simple Lyapunov spectrum 
is still dense, as proved by Arnold and Cong~\cite{ArCon}.
But it ceases to contain an open set: 
they showed that the set of cocyles with one-point spectrum
is also a dense subset.
In this note we show that the set of cocycles with one-point spectrum
is residual, see theorem~\ref{t.generic} below.
That is, generic $L^p$-cocycles, $p<\infty$, have all Lyapunov exponents equal.

The results of Arnold and Cong imply that the Lyapunov exponents 
do not depend continuously on the cocycle.
However, we show that the Lyapunov exponents
have a semicontinuity property.
Using this and Arnold and Cong's density result we obtain
theorem~\ref{t.generic}.

\smallskip

It is interesting to mention that the typical behavior of the spectrum may change drastically
when cocycles more regular than just continuous are considered.
Bonatti and Viana~\cite{BonVi} considered an open class
of H\"older-continuous cocycles over an uniformly hyperbolic base dynamics.
In this setting, they showed that generic cocyles have at least two different 
Lyapunov exponents. 
Even more, the set of cocyles with one-point spectrum has infinite codimension. 
Some of the techniques employed in~\cite{BonVi} originate from 
the theory of products of i.i.d. matrices (e.g.~\cite{GolMar}),
where, as we have mentioned, simple Lyapunov spectrum prevails.

\section{Statement of the results} \label{s.2}

Let $(X,\cF,\mu)$  be a probability space and let $T:X \to X$ an
automorphism preserving the measure $\mu$. 
Let $d \ge 2$ be an integer and denote by $\cG$ the set of all 
(mod $0$ equivalence classes of) measurable maps $A:X\rightarrow \GL$, 
where $\GL$ is endowed with the Borel $\sigma$-algebra. 

We shall consider the set $\GIC$ of all maps $A: X \to \GL$ in $\cG$
satisfying the \emph{integrability condition}:
$$
\int_X \log^+ \|A^{\pm 1}(x)\| \, d\mu(x) < + \infty.
$$
Under this condition, the multiplicative ergodic theorem of Oseledets 
(see~\cite{Ose} or~\cite{Ar})
gives us the Lyapunov exponents (with multiplicity)
$\lambda_1(A,x) \ge \ldots \ge \lambda_d(A,x)$
of the cocycle $A$.

Provided $T$ is ergodic, these functions are constant almost everywhere.
Then we say that $A\in\GIC$ has \emph{one-point Lyapunov spectrum} if all Lyapunov 
exponents are equal, and we denote the set of those cocycles by $\one$.
If the cocycle $A$ takes values in $\SL$ then $A \in \one$ if and
only if all Lyapunov exponents are zero.

\smallskip

Let $1\le p \le \infty$.
Following Arnold and Cong \cite{ArCon},
we endow the set $\cG$ with a ``$L^p$-like'' topology.

Let $\|\mathord{\cdot}\|$ be an operator norm on the set
$\mathrm{M}(d,\R)$ of $d \times d$ matrices.
For any measurable $A : X \to \mathrm{M}(d,\R)$, let 
$$
\|A\|_p = \left(\int_X \|A(x)\|^p \, d\mu(x)\right)^{1/p} 
\quad \text{if $p<\infty$,}
$$
and
$$
\|A\|_\infty = \esssup_{x\in X} \|A(x)\|.
$$
We have $0 \le \|A\|_p \le \infty$.
Next, for $A$, $B\in \cG$, let 
$$
\tau_p(A,B) =
\|A-B\|_p+\|A^{-1}-B^{-1}\|_p .
$$
Then set 
$$
\rho_p(A,B)=
\frac{\tau_p(A,B)}{1+\tau_p(A,B)} \, .
$$
Here it is understood that
$$
\|A-B\|_p = \infty \text{ or } \|A^{-1}-B^{-1}\|_p = \infty
\; \Leftrightarrow \; 
\rho_p(A,B) = 1.
$$
According to \cite{ArCon}, 
$\rho_p$ is a metric on $\cG$ (and therefore on $\GIC$).
Moreover, $(\cG,\rho_p)$ (and hence $(\GIC,\rho_p)$) is complete.

\begin{rem} \label{r.holder}
By the H\"older inequality, we have
$\rho_1(A,B) \le \rho_p(A,B)$ for all $A$, $B \in \cG$, $1\le p \le \infty$.
\end{rem}

\begin{rem} \label{r.GIC}  
If $A\in \GIC$ and $B\in \cG$ with $\rho_p(A,B)<1$ then $B\in \GIC$; see~\cite{ArCon}. 
\end{rem}

Arnold and Cong proved in \cite[theorem~4.5]{ArCon} that 
the set $\one$ is dense in $\GIC$ for any metric $\rho_p$, $1\leq p<\infty$.
We improve this result by showing:
\begin{theorem} \label{t.generic} 
Assume $T$ is ergodic.
Then $\one$ is a residual subset of $\GIC$ in the $L^p$-topology,
for any $1 \le p < \infty$.
\end{theorem}

Let $k=1$,\ldots,$d$.
We will study the following quantities:
$$
\Lambda_k (A) = \int_X 
\left( \lambda_1(A,x) + \cdots + \lambda_k(A,x) \right) \, d\mu(x).
$$
In fact, our main result is:

\begin{theorem} \label{t.semicont}
Let $1 \le p \le \infty$ and endow the set $\GIC$ with the metric $\rho_p$.
Then:
\begin{itemize}
\item[(a)] 
the maps $\Lambda_k : \GIC \to \R$ are upper semicontinuous for all $k=1,\ldots,d$,
that is for every $A\in \GIC$ and $\eps>0$ there exists $0<\delta<1$ such that
if $\rho_p(A,B)<\delta$ then $\Lambda_k(B) < \Lambda_k(A)+\eps$;
\item[(b)]
the function $\Lambda_d : \GIC \to \R$ is continuous.
\end{itemize}
\end{theorem}

From theorem~\ref{t.semicont} and Arnold and Cong's density result,
we shall deduce:

\begin{theorem}\label{t.ptcont}
Assume $T$ is ergodic.
Let $1 \le p < \infty$, and endow the set $\GIC$ with the metric $\rho_p$.
Take $k \in \{1,\ldots, d-1\}$.
Then the map $\Lambda_k : \GIC \to \R$ is continuous at some $A \in \GIC$
if and only if $A \in \one$.
\end{theorem}

Since the set of continuity points of any upper semicontinuous map is residual
(see~\cite{Kur}), theorem~\ref{t.generic} is an immediate corollary of
theorems~\ref{t.semicont} and~\ref{t.ptcont}.

\smallskip

For a description of the continuity points of the functions $\Lambda_k$ 
in the case of bounded or continuous cocycles, see~\cite{BocVi}.

\section{Proofs} \label{s.3}

We are going to uses some basic facts about exterior powers and Lyapunov exponents,
see e.g.~\cite[chapter~3]{Ar}.

\smallskip

Let $k \in \{1, \ldots, d\}$.
We denote by $\wedgek \R^d$ the $k$-th exterior power of $\R^d$.
If $A: \R^d \to \R^d$ is a linear map, then it induces a linear map
$\wedgek A: \wedgek \R^d \to \wedgek \R^d$.
An inner product in $\R^d$ induces an inner product in $\wedgek \R^d$,
and the corresponding operator norms satisfy
\begin{equation}\label{e.norms}
\| \wedgek A \| \le \|A\|^k, 
\quad \forall \, A \in \mathrm{M}(d,\R).
\end{equation}
We fix operator norms as above from now on.

\begin{lemma}\label{l.malandrungen}
Let $A, B \in \GL$ and $k\in\{1,\ldots,d\}$.
Then
$$
\log^+ \|\wedgek B \| \le \log^+ \|\wedgek A \| + k \|B-A\|.
$$
\end{lemma}

\begin{proof}
Using~\eqref{e.norms},
\begin{align*}
\|\wedgek B \|^{1/k} &\le \left(\|\wedgek A \| + \|\wedgek(B-A)\|\right)^{1/k} \\
                     &\le \|\wedgek A \|^{1/k} + \|\wedgek(B-A)\|^{1/k} \\
                     &\le \|\wedgek A \|^{1/k} + \|B-A\|.
\end{align*}
Since $\log^+(x+y) \le \log^+ x + y$ for all $x,y \ge 0$,
the lemma follows.
\end{proof}

We denote, for $A\in \cG$, $x\in X$ and $n\in \N$, 
$$
A^n(x) = A(T^{n-1}x) \cdots A(Tx) A(x).
$$
A basic property of $\Lambda_k$ from which we shall deduce its upper semicontinuity
is the following: 
For all $A\in \GIC$,
\begin{equation}\label{e.inf}
\Lambda_k(A) = \lim_{n\to\infty} \frac{1}{n} \int \log \|\wedgek A^n(x)\|\, d\mu(x) 
             = \inf_n            \frac{1}{n} \int \log \|\wedgek A^n(x)\|\, d\mu(x).
\end{equation}

\smallskip

In what follows, all integrals are meant with respect to the measure~$\mu$.
We shall also need the following measure-theoretic result:

\begin{lemma} \label{l.unif}
Given $f \in L^1(\mu)$ with $f \ge 0$, and $\eta >0$,
there exists $K>0$ such that for all $h$ in $L^1(\mu)$
with $h\ge 0$ and $\|h-f\|_1 < \eta$, we have
$$
\int_{\{h > K\}} h \dmu < 2\eta 
\quad \text{and, consequently,} \quad
\mu(\{h > K\}) < \frac{2\eta}{K}.
$$
\end{lemma}

\begin{proof}
Since $f \ge 0$ is integrable, one can find $\gamma>0$ such that
$$
Z \subset X \text{ measurable, } \mu(Z) < \gamma 
\; \Rightarrow \;
\int_Z f \dmu < \eta. 
$$
Let $K = \gamma^{-1} (\|f\|_1 + \eta)$.
Given $h\in L^1(\mu)$ with $h \ge 0$ and $\|h-f\|_1 < \eta$, 
let $E_h = \{h \le K\}$.
Then
$$
\mu(E_h^c) \le K^{-1} \|h\|_1 < K^{-1}(\|f\|_1 + \eta) = \gamma.
$$
Therefore
$$
\int_{E_h^c} h \dmu \le
\int_{E_h^c} f \dmu + \int_{E_h^c} |h-f| \dmu \le
\eta + \|h-f\|_1 <
2\eta.
$$
This proves the first part of the lemma.
The second part is an immediate consequence.
\end{proof}

Next we prove theorems~\ref{t.semicont} and~\ref{t.ptcont}.

\begin{proof}[Proof of theorem~\ref{t.semicont}]
We first prove part~(a).
By remark~\ref{r.holder}, we may assume $p=1$.

Let $A\in \GIC$ and $\eps>0$ be given.
We will denote 
$$
\hlak(A,x) =
\lambda_1(A,x) + \cdots + \lambda_k(A,x)
$$
Let us first consider the case where the following condition is satisfied:
\begin{equation} \label{e.positive}
\hlak(A,x) \ge 0 \quad \text{for $\mu$-a.e. $x \in X$.}
\end{equation}

From the subadditive ergodic theorem, 
we know that the convergence in 
$$
\hlak(A,x) = \lim_{n \to +\infty} \frac{1}{n} \log \|\wedgek A^n (x) \| ,
$$
takes place almost everywhere and also in $L^1$.
Hence, using~\eqref{e.positive},
$$
\lim_{n \to +\infty} \frac{1}{n} \int_X \log^- \|\wedgek A^n \| \dmu = 0.
$$
We take $N\in\N$ such that (recall~\eqref{e.inf})
$$
\frac{1}{N}\int_X \log^-\|\wedgek A^N \|\dmu < \eps
\quad\text{and}\quad
\frac{1}{N}\int_X \log  \|\wedgek A^N \|\dmu < \Lambda_k(A) + \eps.
$$
Therefore
\begin{equation}\label{e.start}
\frac{1}{N}\int_X \log^+ \|\wedgek A^N \|\dmu < \Lambda_k(A) + 2\eps.
\end{equation}

Let $f = \log^+ \|A\|$ and $\eta = \eps/N$.
Let $K>0$ be given by lemma~\ref{l.unif} applied to $f$ and $\eta$.
Set
\begin{equation} \label{e.delta}
\delta' = \min \left\{ \eta, \eps e^{-K(N-1)}\right\} 
\quad \text{and} \quad
\delta = \frac{\delta'}{1+\delta'} \, .
\end{equation}
Now fix $B \in \cG$ such that $\rho_1(B,A) < \delta$.
Therefore $\|B-A\|_1< \delta'$ and, by remark~\ref{r.GIC}, $B\in\GIC$.
Let $g=\log^+ \|B\|$.
By lemma~\ref{l.malandrungen},
$\|g-f\|_1 \le \|B-A\|_1 < \delta' \le \eta$.
We use lemma~\ref{l.unif} with $h=f$ and $h=g$:
Let
$$
E_f = \{f \le K\} \quad \text{and} \quad E_g = \{g \le K\} \, ; 
$$
then
$$
\int_{E_h^c} h \dmu < 2\eta 
\quad \text{and} \quad
\mu(E_h^c) < \frac{2\eta}{K}
\quad \text{for $h = f$, $g$.}
$$
Set
$$
G = \bigcap_{i=0}^{N-1}T^{-i}(E_f \cap E_g).
$$
Then $G^c$ has small measure:
\begin{equation} \label{e.Gc}
\mu(G^c) \le N \mu(E_f^c \cup E_g^c) 
         <   \frac{4 N \eta}{K}
         =   \frac{4 \eps}{K} \, .
\end{equation}

\smallskip

We are going to bound the expression 
$\frac{1}{N}\int_X \log^+ \|\wedgek B^N\|\dmu$.
To do so, we are going to split the integral in two parts,
$\int_X = \int_{G^c} + \int_G$.
For the first part, we have
$$
\frac{1}{N}\int_{G^c} \log^+ \| \wedgek B^N\| \dmu \le
\frac{1}{N}\sum_{i=0}^{N-1} \int_{T^i(G^c)} \log^+ \| \wedgek B \| \dmu \le
\frac{k}{N}\sum_{i=0}^{N-1} \int_{T^i(G^c)} g \dmu.
$$
For each $i=0$,\ldots,$N-1$ we have, by lemma~\ref{l.unif} and relation~\eqref{e.Gc},
\begin{align*}
\int_{T^i(G^c)} g \dmu 
&=   \int_{E_g^c} g \dmu  
       + \int_{E_g \cap T^i(G^c)} g \dmu \\ 
&<   2\eta + K \mu(E_g \cap T^i(G^c))       
\le 2\eps + K \mu(G^c) < 6\eps \, .
\end{align*}
Hence
\begin{equation}\label{e.1}
\frac{1}{N}\int_{G^c} \log^+ \|\wedgek B^N\| \dmu \le 6k\eps.
\end{equation}

\smallskip

Next we estimate the second part.
Using lemma~\ref{l.malandrungen} and \eqref{e.start} we get 
\begin{equation}\label{e.nomorelogs}
\begin{split}
\frac{1}{N} \int_G\log^+ \|\wedgek B^N \| \dmu 
&\le
\frac{1}{N} \int_G\log^+ \|\wedgek A^N \| \dmu +
\frac{k}{N} \int_G\|B^N-A^N\| \dmu
\\ &\le
\Lambda_k(A) + 2\eps + 
\frac{k}{N}\int_G \|B^N-A^N\| \dmu. 
\end{split}
\end{equation}
To estimate the integral on the right hand side,
we proceed as follows.
Take $x\in G$ and $1 \le i \le N-1$; then
\begin{align*}
\|B^{i+1}(x) - A^{i+1}(x) \| 
&\le \|B(T^i x)\| \, \|B^i(x)-A^i(x)\| + \|B(T^i x) - A(T^i x)\| \, \|A^i(x)\|
\\ 
&\le e^K \|B^i(x)-A^i(x)\| + e^{Ki} \|B(T^i x) - A(T^i x)\|.
\end{align*}
Integrating over $G$ and using $\|B-A\|_1 < \delta'$,
we get
$$
\int_G \|B^{i+1}- A^{i+1} \| \dmu 
\le e^K \int_G \|B^i- A^i \| \dmu 
    + e^{Ki} \delta'.
$$
By induction, we obtain
$$
\int_G \|B^i - A^i \| \dmu \le i e^{K(i-1)} \delta' 
\quad \forall \, i=1, \ldots, N.
$$
Using this relation with $i=N$, \eqref{e.nomorelogs} and~\eqref{e.delta},
we get
\begin{equation}\label{e.2}
\frac{1}{N} \int_G \log^+ \|\wedgek B^N \|\dmu \le
\Lambda_k(A) + (2+k) \eps .
\end{equation}

\smallskip

From~\eqref{e.1} and~\eqref{e.2},
we conclude that
$$
\Lambda_k(B) \le
\frac{1}{N} \int_X \log \|\wedgek B^N \|\dmu \le
\frac{1}{N} \int_X \log^+ \|\wedgek B^N \|\dmu \le
\Lambda_k(A) + (2+7k) \eps.
$$
This proves part~(a) of the theorem in the case condition~\eqref{e.positive}
is satisfied.
(Replace $\eps$ with $\eps/(2+7k)$ along the proof.)

\smallskip

Next we consider the general case.
Again, let $A\in \GIC$ and $\eps>0$.
For $a>0$, consider the $T$-invariant set
$L_a = \{ \hlak(A,x) < -a\}$.
Choose $a$ large enough so that
\begin{equation} \label{e.a}
\int_{L_a} \log^+ \|\wedgek A\| \dmu      < \eps \quad\text{and}\quad
\int_{L_a} \hlak(A,x) \, d\mu(x) > - \eps.
\end{equation}

The cocycle $e^a A$ restricted to $L_a^c$ satisfies condition~\eqref{e.positive}.
If $B \in \cG$ is $\rho_1$-sufficiently close to $A$ then
$\rho_1(e^a A, e^a B)$ will also be small and therefore, by the case already considered, 
$$
\int_{L_a^c} \hlak(e^a B, x) \,d\mu(x) \le
\int_{L_a^c} \hlak(e^a A, x) \,d\mu(x) + \eps,
$$
that is,
$$
\int_{L_a^c} \hlak(B, x) \,d\mu(x) \le
\int_{L_a^c} \hlak(A, x) \,d\mu(x) + \eps.
$$

On the other hand, since $L_a$ is invariant,
$$
\int_{L_a} \hlak(B, x) \,d\mu(x) 
= \inf_n \frac{1}{n} \int_{L_a} \log^+\|\wedgek B^n\| \dmu  
\le \int_{L_a} \log^+\|\wedgek B\| \dmu  .
$$
We may also assume that $\|B-A\|_1<\eps$.
Then, lemma~\ref{l.malandrungen} and~\eqref{e.a},
we have
\begin{multline*}
\int_{L_a} \hlak(B, x) \,d\mu(x) 
\le \int_{L_a} \log^+\|\wedgek A\| \dmu + k \int_{L_a} \|B-A\| \dmu  \\ 
\le (1+k) \eps 
\le \int_{L_a} \hlak(A, x) \,d\mu(x) + (2+k)\eps.
\end{multline*}
We conclude that $\Lambda_k(B) < \Lambda_k(A) + (3+k)\eps$.
This completes the proof of part~(a) of the theorem.

\smallskip

Part~(b) is an easy consequence of part~(a):
The functions
\begin{equation}\label{e.lower}
A \mapsto \widetilde{\Lambda}_k (A) = - \Lambda_k (A^{-1}) =
\int_X \left( \lambda_{d-k+1}(A,x) + \cdots + \lambda_d(A,x) \right) \, d\mu(x)
\end{equation}
are lower-semicontinuous.
In particular, the function $\Lambda_d (\mathord{\cdot}) = \widetilde{\Lambda}_d (\mathord{\cdot})$
is continuous.
\end{proof}

\smallskip

\begin{proof}[Proof of theorem~\ref{t.ptcont}]
Let $A \in \one$ and $\eps>0$.
If $B$ is sufficiently $\rho_p$-close to $A$ then (recall~\eqref{e.lower})
$$
\Lambda_1 (B)             < \Lambda_1 (A) + \eps/k
\quad\text{and}\quad
\widetilde{\Lambda}_1 (A) > \widetilde{\Lambda}_1 (B) - \eps/k.
$$
Therefore
$$
\Lambda_k(B) 
\begin{cases}
\le k \Lambda_1(B)
<   k \Lambda_1(A) + \eps
=   \Lambda_k(A) + \eps, \\
\ge k\widetilde{\Lambda}_1 (B)
>   k\widetilde{\Lambda}_1 (A) - \eps
=   \Lambda_k (A) - \eps.
\end{cases}
$$
showing that $\Lambda_k$ is continuous at~$A$.

\smallskip

Conversely, take $A' \in \GIC$ and assume 
$\Lambda_k$ is continuous at $A'$.
By the theorem of Arnold and Cong~\cite{ArCon}, there exists
a sequence $B_n$ in $\one$ converging to $A'$ in the $\rho_p$-metric.
Since $\Lambda_d$ is also continuous at $A'$, we have
$$
\frac{\Lambda_k(A')}{k} =
\lim_{n \to \infty} \frac{\Lambda_k(B_n)}{k} =
\lim_{n \to \infty} \frac{\Lambda_d(B_n)}{d} =
\frac{\Lambda_d(A')}{d} \, .
$$
This implies $A' \in \one$, because $k<d$.
\end{proof}

As we have mentioned, theorem~\ref{t.generic} follows from
theorems~\ref{t.semicont} and~\ref{t.ptcont}.

\begin{ack}
We would like to thank the referee for his corrections.
\end{ack}


\vfill

{\footnotesize
\noindent Alexander Arbieto ({\tt alexande{\@@}impa.br})
\hfill
Jairo Bochi ({\tt bochi{\@@}impa.br})

\smallskip

\noindent IMPA, Estrada D. Castorina 110, Jardim Bot\^anico, 22460-320 Rio de Janeiro, Brazil
}

\end{document}